\documentclass[letterpaper, 10 pt, conference]{ieeeconf}  % Comment this line out
\IEEEoverridecommandlockouts                              % This command is only
\renewcommand{\baselinestretch}{1}
\usepackage{graphicx}
\usepackage{psfrag}
\usepackage{amsmath}
\usepackage{latexsym}
\usepackage{amssymb}
\usepackage{graphicx}

\def\QED{\begin{flushright} $\Box$ \end{flushright}}

\newcommand{\be}{\end{eqnarray*}}
\newcommand{\ee}{\end{eqnarray*}}
\newcommand{\ben}{\begin{eqnarray}}
\newcommand{\een}{\end{eqnarray}}
\newtheorem{lemma}{Lemma}[section]
\newtheorem{remark}[lemma]{Remark}
\newtheorem{theorem}[lemma]{Theorem}

\begin{document}

\title{\LARGE \bf
A decomposition technique for  pursuit evasion games with many pursuers
}

\author{Adriano Festa and Richard B. Vinter% <-this % stops a space
% <-this % stops a space
\thanks{A. Festa is with Imperial College of London, EEE Department,  London SW7 2AZ, UK.
        {\tt\small a.festa@imperial.ac.uk}}%
\thanks{R. Vinter is with Imperial College of London, EEE Department,  London SW7 2AZ, UK.
        {\tt\small r.vinter@imperial.ac.uk}}%
}

\maketitle
\thispagestyle{empty}
\pagestyle{empty}

%%%%%%%%%%%%%%%%%%%%%%%%%%%%%%%%%%%%%%%%%%%%%%%%%%%%%%%%%%%%%%%%%%%%%%%%%%%%%%%%
\begin{abstract}
Restrictions on memory storage impose ultimate limitations on the dimensionality of differential games problems for which optimal strategies can be computed via direct solution of the associated Hamilton-Jacobi-Isaacs equations. It is of interest therefore to explore whether, for certain specially structured differential games of interest, it is possible to decompose the original problem into a family of simpler, lower dimensional, differential games. In this paper we exhibit a class of single pursuer-multiple evader  games for which a reduction in complexity of this nature is possible. The target set is expressed as a union of smaller, sub-target, sets. The individual differential games in the family are obtained from the original problem by taking as target set an element in the family of sub-target sets, in place of the original target set. We can exploit geometric features of the dynamic constraints and constraints of the problems arising in this way to reformulate them as lower dimensional, simpler to solve, problems. We give conditions under which the value function of the original can be characterized as the lower envelope of the value functions for the simpler problems and how optimal strategies can be constructed from those for the simpler problems. The methodology is illustrated by several examples.
%the simpler problemsThe value function for the original problem is expressed as the lower envelope of the We give conditions under which 

%of the sub-target sets to reduce the dimensionality of the

%by replacing as a union of  sets and  in the decomposition resulti, members of the family of decomposed problemsThe problems considered are pursuer evader games, involving one evader and multiple pursuers. 

%the  for differential games problems of high state dimension Computation of optimal strategies for differential games problems, based on solution of the Data storage constraints limit the dimensionality of differential gamesThe computation place a restriction of the computation of optimal strategies in differential games Here we present a decomposition technique for a class of differential games. The technique consists in a decomposition of the target set  which produces, for geometrical reasons, a decomposition in the dimensionality of the problem. Using some elements of Hamilton-Jacobi equations theory, we find a relation between the  regularity of the solution and the possibility to decompose the problem. We use this technique to solve some cases of a pursuit evasion game with multiple agents.

\end{abstract}

\section{Introduction}
Interest  in Pursuit-Evasion differential games ({\it PE games}) involving several players dates back to the  1960's.  Progress in this field is documented in  the classical differential games literature,  which includes the books by Isaacs \cite{I65}, Pontryagin \cite{P88}, Friedman \cite{F71}, Krasovskii and Subbotin \cite{KS88}.
Constructing optimal strategies for each player, finding the value of the game, deriving optimality conditions for the trajectories and establishing conditions for solvability of the game are typical objectives.\par
The case of several agents, although it can be considered as a special case of the general framework for the problem, has been addressed separately by many authors. An early significant contribution was that of Pshenichnii \cite{P76}, who studied the problem with many pursuers and equal speeds and derived necessary and sufficient conditions for solutions to this problem. Ivanov and Ledyaev \cite{IL81} subsequently studied the optimal pursuit in higher dimensional state spaces, with several pursuers and geometrical constraints. Through the study of an auxiliary problem, related to the interaction between one pursuer and the evader and using a Lyapunov function, they obtained sufficient conditions of optimality. Chodun \cite{C89} and more recently Ibragimov \cite{I05} used the same approach to solve a collection of one to one problems (one pursuer one evader)  for problems with simple dynamic constraints.\par
Differential games, and multi-agent PE games in particular, have found applications in a variety of fields, for example in mathematical economics, where a game has been constructed to model the  relation between agents \cite{J86, E11}, in robotics where typically the emphasis is real time solutions and efficient computational methods. \cite{HW98, VSH02}.\par

Our proposed approach is to decompose the original game into a family of simpler, lower-dimensional games. This is achieve by expressing the original set as a union of smaller, target  subsets. The individual differential games in the family result from the replacing the original target set by each of the target subsets. Special geometric features of the dynamic constraints and constraints can be exploited to reduce the complexity of the new problems generated in this way. Making use of verification techniques originally proposed by Isaacs \cite{EK72,BRP99}, properties of viscosity solutions and techniques of nonsmooth analysis \cite{BCD97,CS04}, we give we a lower-envelope characterization of the value function, and show how to construct optimal strategies for the original problem from those for the  simpler problems. A similar decomposition technique was used in \cite{FKV12}.

%Our proposed approach is different from earlier ones, though it retains consideration of an auxiliary problem. We invoke the verification theorems involving Hamilton-Jacobi-Isaacs equations, originally proposed by Isaacs \cite{EK72,BRP99}, and the theory of viscosity solutions \cite{BCD97,CS04} to validate a decomposition technique. The approach is motivated by a decomposition procedure proposed in \cite{FKV12}, addressing  a different class of differential games. While the problems studied here and in this earlier paper are different in detail, it is possible, a decomposition of target set, in both cases, greatly simplies analysis of the problem.\par
%Initially we investigate the proposed decomposition under a technical hypothesis related to the convexity  of the Hamiltonian on the dual variable. This hypothesis, while somewhat restrictive, is nonetheless satisfied by some PE games of interest, as we show.
%\end{document}
\section{The Hamilton-Jacobi-Isaacs approach to pursuit-evasion games}
The state $y$ of a dynamic system, partitioned as $n$-vector components $y=(y_{1},\ldots,y_{m})$, is governed by the equations
 \begin{equation}\label{DYN}
%\left\{
\begin{array}{l}
y_1'(t)=-g_1(y(t))a_1(t)+h_1(y(t))b(t)+l_1(y(t))
\\
\vdots\\
%y_2'(t)=-g_2(y(t))a_2(t)+h_2(y(t))b_2(t)+l_2(y(t))\\
%\hdots \\
y_{m}'(t)= -g_{m+1}(y(t))a_{m}(t)
\\
\hspace{0.9in} +h_{m}(y(t))b(t)+l_{m}(y(t))\;.
\end{array}
\end{equation}
in which $g_{i}(.),h_{i}(.): \mathbb R^{n} \rightarrow \mathbb R$, $i=1,\ldots,m+1$ and $l_{i}(.):\mathbb R^{n} \rightarrow \mathbb R^{n}$  are given functions.
\vspace{0.05 in}

%for $i=1,\ldots, m+1$.
For each $i$ we interpret the state component $y_{i}$ to be the relative position of the evader with respect to the $i$'th pursuer.
% as the state of the evader, and the remaining  components $y_{1},\ldots,y_{m}$ as the states of the $m$ pursuers.
% and the $n$-vector component $y_{m+1}$ with the single evader. 
The dimension of the state vector $y$  is therefore $N=n \times m$. The $n \times m$ vector $a=(a_{1},\ldots, a_{m})$ comprises the $n$-vector pursuer controls and  $b$ is the $n$-vector evader control. 

The pursuer and evader controls $a$ and $b$ take values in the sets
 \begin{equation}\label{H3}
\begin{array}{l}
A=B_n(0,\rho_{a})\times \ldots \times B_n(0,\rho_{a}) ,\quad
B=B_n(0,\rho_b) 
%\times \ldots \times B_n(0,\rho_b)
\\
%g_i(x)\rho_{a}-h_i(x)\rho_b\geq0, \quad \forall x\in \mathbb{R}^N, \quad \forall i \in\mathcal{I}.
\end{array} 
%\right.
\end{equation}
for some given numbers $\rho_{a},\,\rho_{b} >0$. Define
\begin{eqnarray}
\mathcal{A}:= \{\mbox{ meas. functions } a:[0,+\infty)\rightarrow A \}\\
\mathcal{B}:= \{\mbox{ meas. functions } b:[0,+\infty)\rightarrow B \}
\end{eqnarray}

It is assumed that
\vspace{0.05 in}

\noindent
(H): $g_{i}(.),h_{i}(.), l_{i}(.)$, $i=1,\ldots m+1$ are Lipschitz continuous, and
\vspace{0.05 in}

\noindent
$g_i(x)\rho_{a}-h_i(x)\rho_b -|l_{i}(x)| > 0, \quad \forall x\in \mathbb{R}^N, \quad \forall i 
%\in\mathcal{I}
$.
\vspace{0.1 in}

%\begin{remark}
%In this formulation, information about the allowable velocities for every pursuer is implicit in the 
%contained in 
%the functions specifying 
%the dynamics and the control constraint sets. 
Note that the maximum allowable magnitude of the velocity of the $i$'th pursuer depends on the states of all the pursuers  and the evader, via the functions $g_{i}(.)$, $h_{i}(.)$ and $l_{i}(.)$. We may assume, without loss of generality that the $g_{i}$'s and $h_{i}$'s are positive functions, since the constraint sets $A$ and $B$ are symmetric about the origin. 
%\end{remark}
%where
 %$y_i\in\mathbb{R}^n$, $y=(y_1^t,y_2^t...y_m^t)^t\in\mathbb{R}^N$ (we have that $N=n\cdot m$) and each equation is related to  a %pursuer; evidently $g_i:\mathbb{R}^n\rightarrow\mathbb{R}^+$, $h_i:\mathbb{R}^n\rightarrow\mathbb{R}^+$  and $l_i:\mathbb{R}^n\rightarrow\mathbb{R}^n$.  
 \vspace{0.1 in}
 
 Take the target set to be
\begin{equation}
\label{target}
\mathcal{T}=\{(y_1,y_2,...,y_{m})\in \mathbb{R}^N: \min_{i\in \{1,2...m\}}|y_i|\leq r \}\,,
\end{equation}
in which   
%$|.|$ denotes the  Euclidean norm, for some 
$r \geq 0$ is a specified number.  
%$\mathbb{R}^n$. 
%We also specify certain target set associated with the individual pursuers, namely
%\begin{equation}
%\label{target1}
%\mathcal{T}_i=\{(y_1^t,y_2^t,...,y_{m}^t)^t\in \mathbb{R}^N: |y_i|\leq r \} \quad i\in\{1,2...m\}.
%\end{equation}
Denote by $y_x(.)=y_x(.,a,b)$ the solution of \eqref{DYN}, for given initial state $x$ and controls $a(.)\in {\cal A}$ and $b(.) \in {\cal B}$.  Define the hitting time for (given $x$, $a$ and $b$) to be 
 \begin{equation}
t_x(a,b):=\left\{
\begin{array}{l}
\min\{t:y_x(t;a,b)\in\mathcal{T}\} \\
+\infty  \quad \hbox{if }y_x(t;a,b)\notin \mathcal{T}  \;
%\phantom{aa}
\forall t.
\end{array} \right.
\end{equation}
The player with control $a$ (comprising the pursuers) seeks to minimize the hitting time, while the evader player, with control  $b$, seeks to maximize it.  It is convenient to transform the hitting time cost by means of the mapping
% classical construction  tranequivalent problem results from transforIt is convenient to renormalize apply a classical re-normalization, in which the hitting time cost is transformed according to
%by means of the nonlinear transformation
 \begin{equation}\label{Kruz}
\psi(u):=\left\{
\begin{array}{ll}
1-e^{-u} & \hbox{if }u<+\infty \\
1 & \hbox{if }u=+\infty \;.
\end{array} \right.
\end{equation}
The cost becomes: 
%giving rise to the cost functional
%and consider the discounted cost functional 
\begin{equation}\label{cost}
J(x,a,b)=\psi(t_x(a,b))=\int_0^{t_x}e^{-s}ds.
\end{equation}
The transformation modifies the value function, but leaves unaltered the optimal strategies, owing to the fact that the transformation is monotone.\par
% of the game with cost $J$ we can easily recover the value function for the original one, by simply applying the inverse transformation $\psi^{-1}(v)=-\ln(1-v)$.  
Full specification of the differential game requires precise description of the nature of the control strategies involved. For this purpose, we follow the Elliot/Kalton approach, based on the concept of `non-anticipative' strategies. The set of such strategies for the $a$-player is 
%\begin{multline}
%\Gamma:=\{\alpha:\mathcal{B}\rightarrow \mathcal{A}: t>0, b(s)=\tilde{b}(s) \hbox{ for all }s\leq t\\
 %\hbox{ implies } \alpha[b](s)=\alpha[\tilde{b}](s) \hbox{ for all }s\leq t \}\;.
%\end{multline}
\begin{multline}
\Gamma:=\{\alpha:\mathcal{B}\rightarrow \mathcal{A}: t>0, b(s)=\tilde{b}(s) \hbox{ for all }s\leq t\\
 \hbox{ implies } \alpha[b](s)=\alpha[\tilde{b}](s) \hbox{ for all }s\leq t \}\;.
 \end{multline}
 \begin{multline}
\Delta:=\{\beta:\mathcal{A}\rightarrow \mathcal{B}: t>0, a(s)=\tilde{a}(s) \hbox{ for all }s\leq t\\
 \hbox{ implies } \beta[a](s)=\beta[\tilde{a}](s) \hbox{ for all }s\leq t \}\;.
 \end{multline}
%The set of non anticipating strategies for the $b$-player is defined in an analogous way, and is denoted by $\Delta$. 
The upper and lower values of the game are
\begin{eqnarray*}
&&
u^{+}(x):=\sup_{\beta\in\Delta} \inf_{a\in\mathcal{A}} J(x,a,\beta[a]) 
\\
&&u^{-}(x):=\inf_{\alpha\in\Gamma} \sup_{b\in\mathcal{B}}J(\alpha[b],b)\;.
\end{eqnarray*}
Under the stated hypotheses, and in view of the fact the Isaac's condition is satisfied, the upper and lower values coincide, for arbitrary initial state $x$. The common value defines the value function  $u(.)$ for the  game. Thus $u(x)=u^{+}(x)=u^{-}(x)$ for all $x$.
% Consequently the upper and lower values for the game coincide   \cite{BRP99} and define the value function:
%are, in this case, coincident and they are defined as
%\begin{equation}
%v(x):=\sup_{\beta\in\Delta} \inf_{a\in\mathcal{A}} J(x,a,\beta[a])=\inf_{\alpha\in\Gamma} \sup_{b\in\mathcal{B}}J(\alpha[b],b)\;.
%\end{equation}
%The fact that the lower and upper values are the same follows from the fact that, owing to the special nature of the dynamics and the cost, the Isaacs' condition is satisfied.  For details see \cite{BRP99}.\par
%Let us recall key features of the theory of this class of games. 
\vspace{0.1 in}

{\it It can be shown that, under the state assumptions, the value function is the unique uniformly continuous function $u(.):\mathbb R^{N}\rightarrow \mathbb R$, which vanishes on on ${\cal T}$ and which is a viscosity solution on $\mathbb R^{N} \backslash {\cal T}$ of  the following Hamilton-Jacobi-Isaacs equation
 \begin{equation}\label{P}
\left\{
\begin{array}{ll}
u(x)+H(x,Du(x))=0 & x\in\mathbb{R}^N\setminus \mathcal{T}\\
u(x)=0 & x\in\partial\mathcal{T}
\end{array} \right.
\end{equation}
 where 
%\begin{multline} 
%H(x,p):=\\\
%\max_{a\in A}\min_{b\in B}\left\{\left(g(x)a-h(x)b-l(x)\right)\cdot p \right\}-1\\
%=\min_{b\in B}\max_{a\in A}\left\{\left(g(x)a-h(x)b-l(x)\right)\cdot p \right\}-1.
%\end{multline}
The Hamiltonian is
%for $p_i\in \R^m$
% $H(x, p=(p_{1},\ldots,p_{m} ))$
%\begin{eqnarray*}
%&&
%\hspace{-0.2 in}H(x,p)=\max_{a\in A} \min_{b\in B}\left\{\left( (g_1(x)a_1-h_1(x)b_1-l_1(x) ),\right.\right. \\
%&& \hspace{-0.2 in}\left.\left.(g_2(x)a_2-h_2(x)b_2-l_2(x) ),... \right.\right.\\
%&&\hspace{-0.2 in}
% \left.\left.(g_m(x)a_m-h_m(x)b_m-l_m(x) )\right)\cdot p\right\}-1
%\\
%\end{multline}
%\begin{multline}\label{pp1}
%&&\hspace{-0.2 in}
%\quad=\,\sum_{i=1}^m g_i(x)\max_{a_i\in B_n(0,\rho_a)}\left\{(a_i)\cdot p_i\right\} \\
%&&\hspace{-0.2 in}
%\quad-\max_{b\in B}\left\{\sum_{i=1}^m h_i(x)\left(b_i\right)\cdot p_i \right\}-\sum_{i=1}^m \left(l_i(x)\right)^t\cdot p_i -1.\hspace{-0.2 in}
%\end{eqnarray*}
\begin{eqnarray*}
&&
\hspace{-0.2 in}
H(x,p=(p_{1},\ldots,p_{m}))\;=
\\
%\max_{a\in A} \min_{b\in B}\left\{\left( (g_1(x)a_1-h_1(x)b_1-l_1(x) ),\right.\right. \\
%&& \hspace{-0.2 in}\left.\left.(g_2(x)a_2-h_2(x)b_2-l_2(x) ),... \right.\right.\\
%&&\hspace{-0.2 in}
% \left.\left.(g_m(x)a_m-h_m(x)b_m-l_m(x) )\right)\cdot p\right\}-1
%\\
%\end{multline}
%\begin{multline}\label{pp1}
%&&\hspace{-0.2 in}
&&\, \rho_{a} \sum_{i=1} g_i(x)|p_{i}| - \rho_{b} |\sum_{i=1} h_i(x) p_{i}| \\
&& \hspace{1.0 in}
- \sum_{i=1} l_i(x)\cdot p_i -1\;.
%
%\max_{a_i\in B_n(0,\rho_a)}\left\{(a_i)\cdot p_i\right\} \\
%&&\hspace{-0.2 in}
%\quad-\max_{b\in B}\left\{\sum_{i=1}^m h_i(x)\left(b_i\right)\cdot p_i \right\}-\sum_{i=1}^m \left(l_i(x)\right)^t\cdot p_i -1.\hspace{-0.2 in}
\end{eqnarray*}

Furthermore, the value function is locally Lipschitz continuous.}
\vspace{0.05 in}

\noindent
Solving this equation for the value function of the game, yields the optimal strategy of each player for initial state $x_0$ as $a(t)=S(y_{x_0}(t))$ and $b(t)=W(y_{x_0}(t))$ where
\vspace{0.05 in}

\noindent
$S(z)\in \underset{a\in A}{\operatorname{argmax}}\, \underset{b\in B}\min\{(g(x)a-h(x)b-l(x))\cdot Dv(x)\}$
\vspace{0.05 in}

\noindent
$W(z)\in \underset{b\in B}{\operatorname{argmin}}\,\underset{a\in A} \max\{(g(x)a-h(x)b-l(x))\cdot Dv(x)\}.$
\vspace{0.1 in}

\noindent{\bf Example 1.}  We examine a simple, preliminary example. Consider the pursuit-evasion game with two pursuers $p_1, p_2$ and one evader $e$, and the positions of each evolve in 1D space.
% with various velocities. 
We take the dynamics to be
 \begin{equation}
\left\{
\begin{array}{l}
p_1'=\frac{2}{3}a_1\\
p_2'=a_2\\
e'=\frac{b}{2}\\
p_1(0)=p_1^0\\
p_2(0)=p_2^0\\
e(0)=e^0
\end{array} \right.
\end{equation}
where $a_1,a_2,b\in B(0,1)=[-1,1]$, $p_1,p_1,e\in\mathbb{R}$. The speed constraint on pursuers is assumed to be greater than that of the evader; this  ensures that the optimal hitting time is finite for an arbitrary initial finite state.
%\vspace{0.05 in}

Let us now consider a reduced formulation. We translate the origin to the position of the e-player, so the $i$'th state variable becomes  the position of the evader relative to that of the $i$'th pursuer. There results the reduced dynamics
 \begin{equation}\label{dyn3}
\left\{
\begin{array}{l}
y_1'=-\frac{2}{3}a_1+\frac{b}{2}\\
y_2'=-a_2+\frac{b}{2}\\
y_1(0)=p_1^0-e^0\\
y_2(0)=p_2^0-e^0
\end{array} \right.
\end{equation}
here we have $a_1,a_2,b\in B(0,1)$, $y_1,y_2\in (-\infty,+\infty)$.  The target set $\mathcal{T}$ becomes the union of neighbourhoods the two axes:
%in this case is the neighborhood of the two axes
%\begin{equation}
%solution, but  is  not necessary 
%for this purpose,.
%request for a well position of the problem, 
%in view of  transformation \eqref{Kruz} that is introduced. 
capture occurs when $\min_{i\in\{1,2\}}|p_i-e|\leq r$, for some specified $r\geq 0$. 
%Note that, in this problem, the state space has dimension two. \par
\begin{figure}[ht]
\begin{center}
\includegraphics[height=1.6cm]{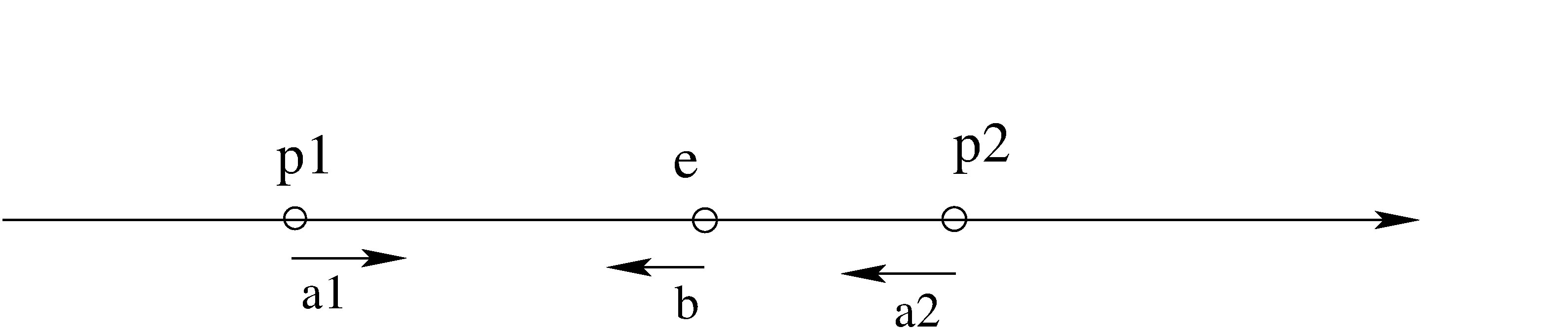}
\caption{A monodimensional PE game with two pursuers.}
\end{center}
\end{figure}

%\begin{equation}
%\mathcal{T}:=\{(x_1,x_2)\in\mathbb{R}^2: \min_{i\in\{1,2\}}|x_i|\leq r \}.
%\end{equation}
The HJI equation  is
 \begin{equation}\label{HJIex11}
\left\{
\begin{array}{l}
u(x)+\max\limits_{a_1,a_2}\min\limits_{b}\left\{(\frac{2}{3}a_1-\frac{b}{2},a_2-\frac{b}{2})\cdot Du(x)\right\}=1 \\
\hspace{1.4 in} x\in[0,+\infty]^2 \setminus\mathcal{T}\\
u(x)=0 \hspace{0.9 in}

x\in\partial\mathcal{T}
\end{array} \right.
\end{equation}

\par
Observe that, with this formulation, the dimension of the domain in which we seek to solve \eqref{HJIex11} is two. The rapid growth in complexity of the problem with  increase in the number of the pursuers places severe restrictions on the computability of this equation for state dimensions greater than three. We now propose a decomposition technique to address this challenge.\par

\section{A decomposition technique for the $m$ pursuer problem}
In this section we present a decomposition technique to overcome the complexity of the multiple pursuer game introduced in Section 2,  for high state dimensions. The key idea is to decompose the target set ${\cal T}$ as a union of smaller sets ${\cal T}_{i}$, $i=1,\ldots, m$:
%This will require the imposition of  some additional hypotheses on the dynamics and on the sets $A$ and $B$. \par
%Assume the following dynamics, which may be viewed as a special case of  \eqref{DYN} and \eqref{dyn2}, as remarked in example 1.
 %\begin{equation}\label{dyn}
%\left\{
%\begin{array}{l}
%y_1'(t)=-g_1(y(t))a_1(t)+h_1(y(t))b_1(t)+l_1(y(t))\\
%y_2'(t)=-g_2(y(t))a_2(t)+h_2(y(t))b_2(t)+l_2(y(t))\\
%\hdots \\
%y_m'(t)=-g_m(y(t))a_m(t)+h_m(y(t))b_m(t)+l_m(y(t))
%\end{array} \right.
%\end{equation}
%where
 %$y_i\in\mathbb{R}^n$, $y=(y_1^t,y_2^t...y_m^t)^t\in\mathbb{R}^N$ (we have that $N=n\cdot m$) and each equation is related to  a pursuer; evidently $g_i:\mathbb{R}^n\rightarrow\mathbb{R}^+$, $h_i:\mathbb{R}^n\rightarrow\mathbb{R}^+$  and $l_i:\mathbb{R}^n\rightarrow\mathbb{R}^n$.  Take the target set to be
%\begin{equation}
%\label{target}
%\mathcal{T}=\{(y_1^t,y_2^t,...,y_{m}^t)^t\in \mathbb{R}^N: \min_{i\in \{1,2...m\}}|y_i|\leq r \}
%\end{equation}
%where  $|\cdot|$ denotes the  Euclidean norm 
%$\mathbb{R}^n$. 
%We also specify the target set for the $i$'th pursuer
\begin{equation}
\label{target1}
\mathcal{T}_i=\{(y_1,y_2,...,y_{m})\in \mathbb{R}^N: |y_i|\leq r \} ,\; i=1,\ldots,m.
\end{equation}
and to relate the value function $u(.)$ for the original game to the value functions $u_{i}(.)$ for the games in which the ${\cal T}_{i}$'s replace ${\cal T}$. The $u_{i}(.)$'s are viscosity solutions to the HJI equation above, with modified boundary condition:
% are viscosity solutions to 
%\vspace{0.05 in}
%
%Consider now the following class of problems, with $i\in\mathcal{I}:=\{1,...m\}\subset\mathbb{N}$,
\begin{equation}\label{P1}
\left\{
\begin{array}{ll}
u_i(x)+H\left(x,Du_i(x)\right)=0& x\in\mathbb{R}^N\setminus\mathcal{T}_i\\
u_i(x)=0 & x\in\mathcal{T}_i
\end{array} \right.
\end{equation}

Under the hypotheses, the $u_{i}(.)$'s are Lipschitz continuous.
\vspace{0.05 in}

Define the index set
$$
I(x)\,:=\, \{i\in \{1,\ldots,m\}: u_{i}(x)=\min_i u(i) \}\;.
$$
\begin{theorem} 
\label{t:1}
Assume condition (C) is satisfied:
\vspace{0.05 in}

(C): for arbitrary $x \in \mathbb R^{N}\backslash {\cal T}$, any convex combination 
$\{\lambda_{i}\,|\,  i \in I(x) \}$ and any collection of vectors $\{\xi_{i} \in \partial^{L}u_{i}(x)\,|\,  i \in I(x)   \}$ we have
\begin{equation}
\label{convex}
H(x,\sum_{i}\lambda_{i}\xi_{i})\,\leq\, \sum_{i \in I(x)} \lambda_{i} H(x,\xi_{i})\,.
\end{equation}
\vspace{0.05 in}
Then
$$
u(x)\;=\;  \underset{i}\min\{u_{1}(x),\ldots, u_{m}(x)\} \,
$$
for all  $x \in \mathbb R^{N}\backslash {\cal T}$.
 
\end{theorem}
\vspace{0.08 in}
\begin{remark}
Here the limiting superdifferential  $\partial^{L}u_{i}(x)$ is the set
$$
\partial^{L}u_{i}(x)\; :=\; 
\underset{x' \rightarrow x}{\mbox{lim sup}}\; \partial^{F}u_{i}(x) \;,
$$
in which $\partial^{F}u_{i}(x)$ is the super Frechet differential
\vspace{0.05 in}

$ \quad \partial^{F}u_{i}(x)\,=\,$
\vspace{0.1 in}

$\{
p\in \mathbb{R}^N : \underset{x' \rightarrow x}{\mbox{lim sup}}\,
%_{\substack{y\rightarrow x} 
%\\y\in B(x_0,\rho)}
 \frac{u_{i}(x')-u_{i}(x)-p\cdot  (x'-x)}{|x'-x|}\leq 0
\}\,.$
\vspace{0.1 in}

\noindent
Condition (C) is automatically is satisfied if $H(x,.)$ is a convex function, but is in fact significantly weaker. 
\end{remark}

%However the greater flexibility of condition (C)
%The cost functional is as in \eqref{cost}.
%\par
%We assume
% \begin{equation}\label{H3}
%\left\{
%\begin{array}{l}
%A=B_n(0,\rho_{a})\times B_n(0,\rho_{a})\times ... \times B_n(0,\rho_{a}) ,\\
%B=[B_n(0,\rho_b)]^m,\\
%g_i(x)\rho_{a}-h_i(x)\rho_b\geq0, \quad \forall x\in \mathbb{R}^N, \quad \forall i \in\mathcal{I}.
%\end{array} \right.
%\end{equation}
%where with $A$ is the space  $\{(a_1^t,a_2^t,...,a_m^t)^t\in \mathbb{R}^N: a_i\in\mathbb{R}^n \hbox{ and } |a_i|=\rho_a,\forall i\in\mathcal{I} \}$. $[B_n(0,\rho_b)]^m$ is the set $\{(b_1^t,b_2^t,...,b_m^t)^t\in \mathbb{R}^N: b_i\in\mathbb{R}^n,  b_i=b_j \hbox{ and } |b_i|=\rho_b,\forall i,j\in\mathcal{I} \}$.\par
%\vspace{0.1 in}
%As a first step to establishing such a relation, we note the following property of the Hamiltonian:

%Our goal now is to provide a convenient decomposition of the value function. A key step is summarized as the following lemma: 
%\begin{lemma}\label{convex}
\noindent 
\vspace{0.1 in}
%The proof makes use of the concept of semi-concavity, defined as follows
%Before outlining the proof we remind the basic definition
%\begin{definition}
%We say that $v:\mathbb{R}^N\rightarrow \mathbb{R}$ is \emph{semiconcave} on the open convex set $\omega$ if there exists a constant $C\geq 0$ such that
%\begin{equation}\label{sc}
%\lambda v(x)+(1-\lambda) v(y)\leq v\left(\lambda x + (1-\lambda )y\right)+\frac{1}{2} C \lambda (1-\lambda) |x-y|^2
%\end{equation}
%for all $x,y\in\mathbb{R}^N $ and $\lambda\in[0,1]$.
%\end{definition}
\vspace{0.1 in}

\noindent
%{\it Outline of Proof of Thm. \ref{t:1}}:

%\proof 
\noindent
{\it Proof Outline.}
Define
\begin{equation}
\label{envelope}
\overline{u}(x):=\min\{u_{i}(x);  i\in 1 \ldots m\}\quad \mbox{for all }x 
\end{equation}
Our aim is to show that $\overline{u}(.)$ coincides with $u(.)$.
%,given by (\ref{envelope}),  
Since $\overline{u}(.)$ vanishes on ${\cal T}$, and in view of the facts that $\overline{u}(.)$ is Lipschitz continuous and that $u(.)$ is the unique uniformly continuous viscosity solution to the HJI equation satisfying this boundary condition, it suffices to show that $\overline{u}(.)$ is such a solution.
\vspace{0.05 in}

That $\overline{u}(.)$ is a super (viscosity) solution follows directly from the definition of super solution and the fact that $\overline{u}(.)$ is the lower envelope of a finite number of super solutions. It remains to show then that $\overline{u}(.)$ is a sub solution.  Take any $x \in \mathbb R^{N} \backslash {\cal T}$ and any $\xi \in \partial^{F}(x)$. According to a well-known characterization of sub solutions, it suffices to demonstrate that 
%We must show that
\begin{equation}
\label{valid}
\overline{u}+H(x, \xi) \,\leq\, 0\;.
\end{equation}
But, by the `max rule' for limiting subdifferentials of Lipschitz continuous functions, applied to $-\overline{u}(x)= \underset{i}{\max}(-u)(x)$, we know that, for some convex combination 
%\end{document}
$\{ \lambda_{i}\,|\, i \in I(x) \}    $ 
and set of vectors $\{ \xi_{i} \in \partial^{L}u_i(x)\,|\, i \in I(x)\}$. 
$$
\xi= \underset{i\in I(x)}{\sum} \lambda_{i} \xi_{i}\;.
$$
For each $i \in I(x)$,  there exists sequences $x^{i}_{j}\rightarrow x$ and $\xi^{i}_{j} \rightarrow \xi_{i}$ such that, for each $i$, $\xi\in \partial^{F}(x^{i}_{j})$ as $j \rightarrow \infty$. Since $u_{i}(.)$ is a subsolution, $H(x^{i}_{j},\xi^{i}_{j})\leq 0$. But then, since $H(.,.)$ is continuous,
$$
H(x,\xi_{i})\leq \underset{j}{\mbox{lim sup}}\, H(x^{i}_{j},\xi^{i}_{j} )\,\leq 0\,
$$
It follows from condition (\ref{convex}) that
$$
H(x, \xi)\,=\, H(x, \underset{i}{\sum} \lambda_{i} \xi_{i} )\,\leq\,  \underset{i}{\sum}\lambda_{i} H(x,\xi_{i} ) \;.
$$
then   
\begin{multline}\label{ppp}
\overline{u}+H\left(x, \sum_i \lambda_i p_i\right) \leq \sum_i \lambda_i \overline{u}+\sum_i \lambda_i H(x,p_i) \\
\leq \sum_i \lambda_i \left(u_i+ H(x,p_i)\right) \leq  0.
  \end{multline}
We have confirmed (\ref{valid}). \hspace{3.5cm}{$\square$}

\bigskip

We now discuss the role the decomposition technique, summarized as Theorem \ref{t:1},  in reducing the complexity of the differential game.  \par
Consider again Example 1. The target can be expressed as a union of two sets:  $\mathcal{T}:=\mathcal{T}_1\cup \mathcal{T}_2$, in which $\mathcal{T}_i:=\{(x_1,x_2)\in\mathbb{R}^2 : |x_i|\leq r \}$. For the target $\mathcal{T}_1$ (respectively $\mathcal{T}_2$), the value function is clearly independent of $x_{2}$ (respectively $x_{1}$). So
%second state component,  the second equation of the dynamic \eqref{dyn3} will not affect the value function because we cannot  reach the target using it, and this mean that the value function of the decomposed game with $\mathcal{T}_1$ as target will be constant along $x_2$. As consequence of it, 
$\frac{\partial}{\partial x_2} u_1(x_1,x_2)=0$  and $\frac{\partial}{\partial x_1} u_2(x_1,x_2)=0$. 
$u_{1}(.)$  therefore satisfies
%$$; so we have a simplification of the equation. To get the value function of the decomposed problem, we solve the equation
\begin{figure}[t]
\begin{center}
\includegraphics[height=6cm]{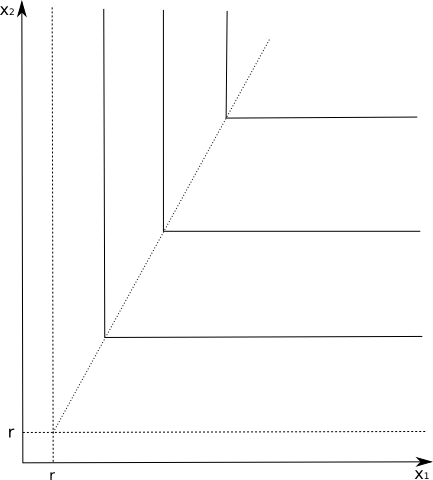}
\caption{Some level sets for the value function of the monodimensional problem with two pursuers}
\end{center}
\end{figure}
\begin{equation}
\left\{
\begin{array}{l}
u_1+\max\limits_{a_1}\min\limits_{b}\left\{-(-\frac{2}{3} a_1+\frac{b}{2})\cdot  \frac{\partial}{\partial x_1}u_1\right\}=1 \\
 \phantom{ffvvvtttttfffff}  
 x_1\in (r,+\infty], x_2\in [0,+\infty)\\
u_1=0 
\phantom{fvvfffff} 
x_1\in [0,r],x_2\in [0,+\infty)
\end{array} \right.
\end{equation} 
where $a_1\in B_n(0,\rho_{a})$ and $b\in B_n(0,\rho_b)$. This is a 1D equation for a fixed $\overline{x}_2$ and is constant for a fixed $\overline{x}_1$. It has solution
%So this is enough simple to solve. We get
\begin{equation}
u_1(x)=1- e^{-6(x_1-r)} \hbox{,  } \psi^{-1}(u_1)=6(x_1-r)
\end{equation}
Similarly
%similar situation is with the second decomposition problem. If we consider just $\mathcal{T}_2$ the first equation of the dynamics does not affect the decomposed problem, so in the same way
\begin{equation}
u_2(x)=1- e^{-2(x_2-r)} \hbox{,  } \psi^{-1}(u_2)=2(x_2-r)\;.
\end{equation}
%Now apply Theorem \ref{t:1}. 

%The Hypotheses of the Thm.  \ref{t:1} are satisfied in this case, in which $g(x)\equiv \left(\frac{2}{3},1\right)$, $h(x)\equiv \left(\frac{1}{2},\frac{1}{2}\right)$ and $\rho_a=\rho_2=1$. This means that $\left(g_1(x)\rho_a-h_1(x)\rho_b,g_2(x)\rho_a-h_2(x)\rho_b\right)=\left(\frac{2}{3}-\frac{1}{2},1-\frac{1}{2}\right)=(\frac{1}{6},\frac{1}{2})$. The Hamiltonian $H(x,.)$ is not convex in the second variable. Noting however that, for given of $x$, gradients of $u_{1}(.)$ and of $u_{2}(.)$ 
%%\end{document}
\noindent
The gradients of these two functions at $x$ are of the form $(k_{1}(x),0)$, $(0, k_{2}(x))$ for non-negative functions $k_{1}(.)$ and $k_{2}(.)$. It is simple to check that condition (C) is satisfied. According to Thm. \ref{t:1},  the value function as lower envelope of $u_{1}(.)$ and $u_{2}(.)$, thus
%
%
%In this way we can build the correct value function of the original problem
\begin{equation}\label{HJIex}
\psi^{-1}(\overline{u})=
\left\{
\begin{array}{ll}
-6(x_1-r)& \hbox{if }x_1\leq \frac{1}{3}x_2+\frac{2}{3}r\\
-2(x_2-r) & \hbox{if }x_1 > \frac{1}{3}x_2+\frac{2}{3}r.
\end{array} \right.
\end{equation} 
%Solving two simpler problems (the dimension of the decomposed problem is lower than the original one because of the geometry of the domain) we get the solution of the problem.\par

\begin{remark}
These optimal strategies are consistent with  intuition. When one pursuer is very close to the evader, that pursuer's location alone affects the strategy of the evader, as expected.  This feature of the solution is evident from the formulae for the value function, which reveal that, for states far from the $x_{1}$ axis  but close to the $x_{2}$ (for example), the value function  coincides with the function $u_i(x)$. 
%More interesting is the fact that also in the connection points between two different decomposed solution,  when two pursuers are equally close to the capture, the value function of the game is well reconstructed.
\end{remark}

The preceding analysis can be generalized to cover a unidimensional problem with $m$ pursuers.
%With a similar reasoning, we can generalize this procedure
%\end{document}
\begin{theorem}\label{decom}
Assume hypotheses (H). 
%\eqref{H3}, 
Denote by  $v_i(.):\mathbb{R}\rightarrow \mathbb{R}$ the solution of the following equation
\begin{equation*}
\left\{
\begin{array}{l}
v_i(x_{i})+\max\limits_{a_i}\min\limits_{b}\left\{f_i(x_{i},a_i,b)\cdot  D v_i(x_{i})\right\}=1 \\
\phantom{tttppppppppppppppppppppttttttttttt} x_{i}\in (r,+\infty]\\
v(x_{i})=0 \phantom{gggggggggggggggggggggg} x_{i}\in [0,r]
\end{array} \right.
\end{equation*} 
where $f_i(x_{i},a_i,b):=g_i(x)a_i-h_i(x_{i})b-l(x_{i})$, $a_i\in B_1(0,\rho_{a})$, $b\in B_1(0,\rho_{b})$.
% and defined the function $U_i(x):\mathbb{R}^N\rightarrow\mathbb{R}$ as
\begin{equation*}
u_i(x=(x_1,x_2...x_i,...x_m))=v_i(x_i) 
%\quad \forall (x_1,x_2,...x_m)\in\mathbb{R}^{N}
\end{equation*}
%we claim that the value function of the original problem is
%$$ u(x)=\min_{i\in\{1,2...m\}}U_i(x).$$
Then the value function $u(.)$ is
$$
u(x) \;=\; \min\{ u_{1}(x),\ldots, u_{m}(x)        \}
$$
\end{theorem}
\vspace{0.1 in}

\proof
It is straightforward to confirm that $u_i(x)$   is a viscosity solution of the problem \eqref{P1}, 
using the fact that $\frac{\partial u_i}{\partial x_j}(x)=0$ for all $i\neq j$.
%, the equation \eqref{P1} reduce to
%$$ U_i(x)+\max\limits_{a_i\in B_1(0,\rho_{a})}\min\limits_{b_i\in B_1(\rho_{b_i})}\left\{f_i(x,a_i,b)\cdot  D U_i(x)\right\}=1 $$
%which is verified for definition by $U_i$. 
Then, in view of the uniqueness of viscosity solutions for \eqref{P1} (\cite{BCD97} Theorem 3.1), we know that $u_i(.)$ is the unique viscosity solution of the decomposed problem $\eqref{P1}$.\par
%Now, we want to use Theorem \ref{t:1} to validate the decomposition procedure. 
%Therefore we we need to prove that condition (C) is always verified. 

Making use of the fact that that $\xi_i\in\partial^L u_i(x)$ has the structure $(0,\ldots,0, k_i ,0,\ldots,0)$, where $k_i\geq 0$, we deduce
Then, for any $\xi_i \in\partial^L u_i(x)$ and $\xi_j \in\partial^L u_j(x)$ with $i,j\in I(x)$, we have
\begin{multline}
H(x,\xi_i+\xi_j)=\\
\rho_a \left( g_i |k_i|+g_j|k_j|\right) -\rho_b |h_i k_i + h_j k_j|-l_i \cdot (0...k_i ...k_j...0)\\
= \left(\rho_a  g_i |k_i| -\rho_b |h_i k_i |-l_i \cdot (0...k_i ...0)\right) \\
+  \left(\rho_a  g_j |k_j| -\rho_b |h_j k_j |-l_j \cdot (0...k_j ...0)\right) \\
= H(x,\xi_i)+H(x,\xi_j).
\end{multline}
this confirms condition (C) of Thm. \ref{t:1} is satisfied. The stated characterization of the value function now follows from Thm. \ref{t:1}.
%\vspace{0.5 in}
%
%\noindent
%It is sufficient to remark that for construction, a $\xi_i\in\partial^L u_i(x)$ has the structure $(0...0\; k_i \;0...0)$, where $k_i\geq 0$. This is due to the monotone nature of the unidimensional minimum time problem, with a connected target. Then, said $\xi_i \in\partial^L u_i(x)$ and $\xi_j \in\partial^L u_j(x)$ with $i,j\in I(x)$, we have
%\begin{multline}
%H(x,\xi_i+\xi_j)=\\
%\rho_a \left( g_i |k_i|+g_j|k_j|\right) -\rho_b |h_i k_i + h_j k_j|-l_i \cdot (0...k_i ...k_j...0)\\
%= \left(\rho_a  g_i |k_i| -\rho_b |h_i k_i |-l_i \cdot (0...k_i ...0)\right) \\
%+  \left(\rho_a  g_j |k_j| -\rho_b |h_j k_j |-l_j \cdot (0...k_j ...0)\right) \\
%= H(x,\xi_i)+H(x,\xi_j).
%\end{multline}
%So, using Theorem \ref{t:1} we get the thesis.
%
\endproof

\section{Examples and numerical tests. }
In this section we solve some higher dimensional problems using the proposed decomposition technique. We note that memory storage constraints impose fundamental limits on the state dimensions of problems for which solutions to the associated HJI equations can be computed directly. Indeed, MATLAB implementations using  a heap-based Java VM system are not feasible for $N$ dimensional arrays, for $N >5$ . This provides the motivation for studying decomposition techniques which are, in principle, applicable for very high dimensional problems.
% solvedessentially limit the dimension
%
 %introduced in the previous section. We will like to focus on the fact that  for the practical point of view, an important bottleneck is in the memory storage.  We had to stop our tests managing 5D arrays because the not so efficient management of the memory (we used as developing program MATLAB which uses  a heap-based system of the Java VM). This problem can be overcome with a more specialist storage of the information. We considered that this point exceeds of the aim of this paper.

\bigskip
{\bf Test 1}\par
\medskip
In this example we study
 pursuit-evasion game on a plane involving several pursuers and one evader, with the help of Thm. 3.1.\par
Consider the problem in a reduced space where the variable $x_i$ is the distance between the evader and the $i$'th pursuer and  every pursuer  has constant maximum speed $g_i(x)\equiv 1$. The evader has constant maximum speed $\displaystyle h_i(x)= 0.9 $. \par

The $i$'th unidimensional problem is, for control constraint sets $a_i\in B_1(0,1)$ and $b\in  B_1(0,1)$, has associated HJI equation
\begin{equation}
\left\{
\begin{array}{l}
v_i(x_i)+\min\limits_{a_i}\max\limits_{b}\left\{(g_i a_i- h_i(x) b)\cdot  D_i v_i(x_i)\right\}=1 \\
\phantom{gggggggggggggggggggggggggg} 
x_i\in (r,+\infty]\\
v(x_i)=0 
\phantom{gggggggggggggggggg} 
x_i\in [0,r].
\end{array} \right.
\end{equation} 
We extend now the $v_i$ solution in the space of dimension 5: $V_i(x):\mathbb{R}^{5}\rightarrow R$ is defined as 
\begin{equation}
u_i(x=(x_1,...x_i,...x_5))=v_{i}(x_i) \;.
%\quad \forall (x_1,...x_5)\in\mathbb{R}^{5}
\end{equation}
The value function of the problem is simply 
\begin{equation}
u(x)=\min\{u_{1}(x),\ldots, u_{5}(x)\}.
\end{equation}

\begin{figure}[th]
\begin{center}
\includegraphics[height=7cm]{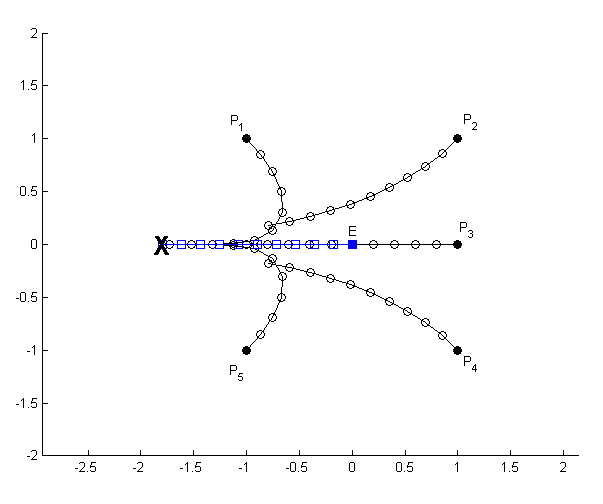}
\caption{Test 1: optimal trajectories at various times. A $X$ indicates the point of capture.} \label{fig4}
\end{center}
\end{figure}

\begin{figure}[th]
\begin{center}
\includegraphics[height=7cm]{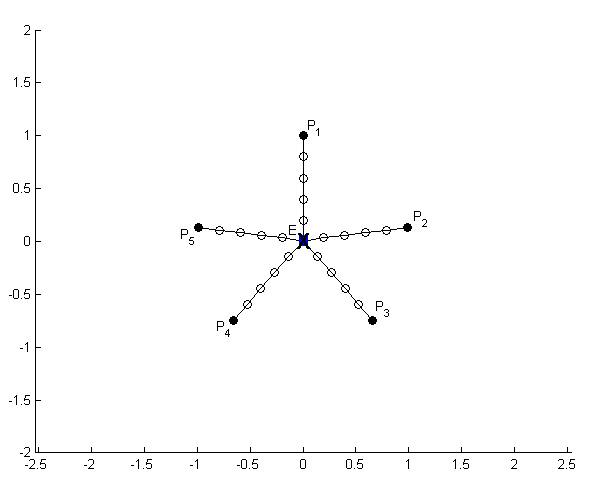}
\caption{Test 1: optimal trajectories at various times. A $X$ indicates the point of capture. The optimal trajectory of the evader is waiting the capture.} \label{fig6}
\end{center}
\end{figure}

Figs.  \ref{fig4}, \ref{fig6} we show simulations for two different with different starting points.

%
%\begin{figure}[th]
%\begin{center}
%\includegraphics[height=7cm]{test2.png}
%\caption{Test1: optimal trajectories at various times. The capture does not happen before the evader reaches the boundary of the  set $[-2, 2]\times[-2, 2]$.} \label{fig5}
%\end{center}
%\end{figure}

\bigskip
{\bf Test 2}\par
\medskip
%In this case we want to consider a different kind of PE game, in particular with a different target set. We can see that the result \ref{t:1} is still effective and the decomposition technique above can be easily adapted to this situations. \par 

\begin{figure}[t]
\includegraphics[height=6.3cm]{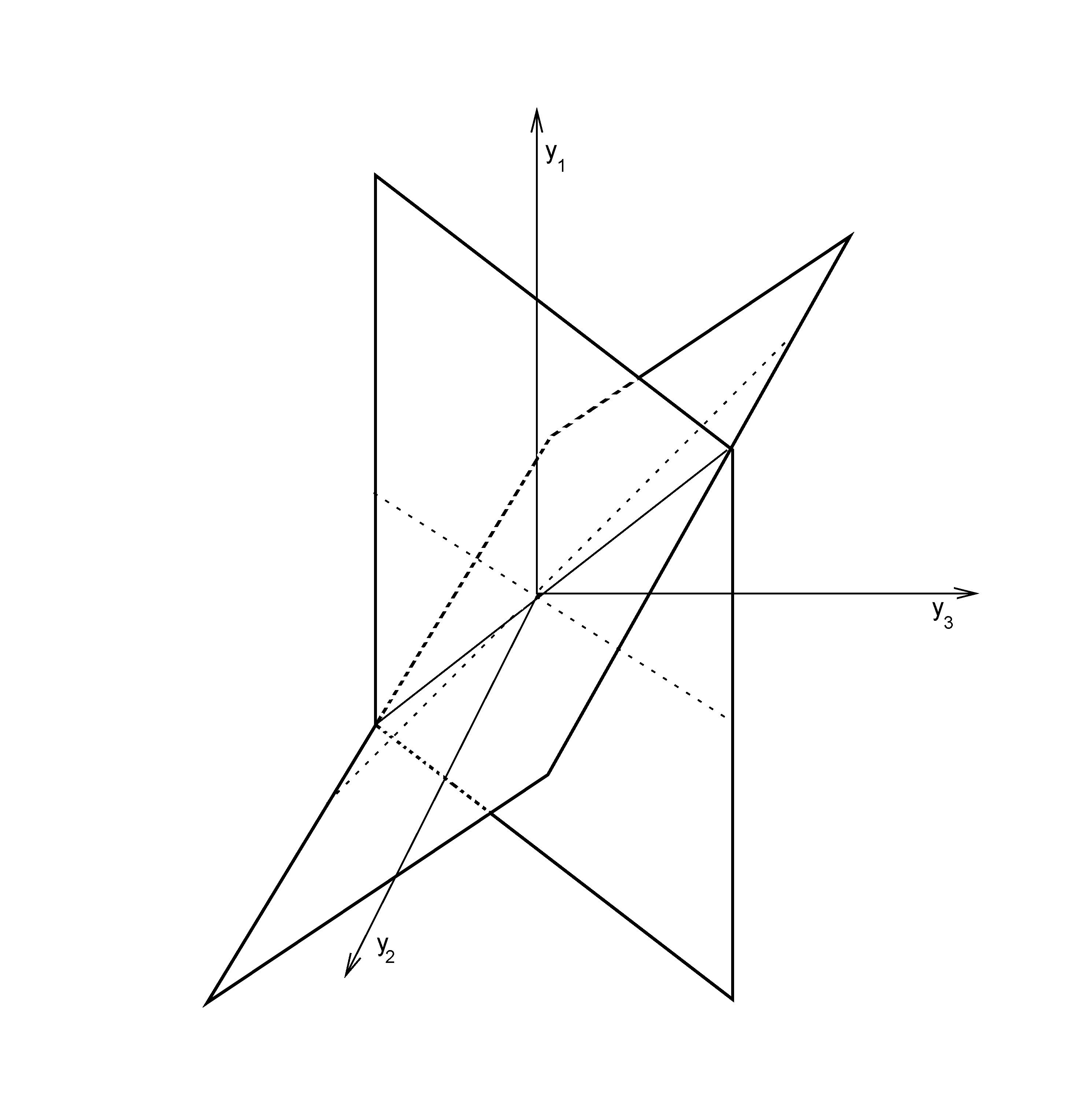}
\caption{Test 2: Target set for $m=2$, $n=1$, $r=0$.} \label{targset}
\end{figure}

In this example, we consider a PE game involving several pursuers and one evader, but now in presence of an obstacle (a river, a different kind of medium) which affects the velocities of the pursuers and the evader. In this case we cannot pass to the reduced coordinates, since this will result in a loss of the information about the position of the agents. Accordingly, the first $m$ block components of the state are associated with the pursuers and the $(m+1)$'th block component with the evader. The positions of the agents, evolving in $\mathbb{R}^n$,  are governed by the following equations
%dynamicsSo, if the space where the agents are moving is $\mathbb{R}^n$ we have the following dynamics, where in the component $m+1$ we will settle the evader:

 \begin{equation}\label{dynex3}
\left\{
\begin{array}{l}
y_1'(t)=-g(y(t))a_1(t)\\
y_2'(t)=-g(y(t))a_2(t)\\
\hdots \\
y_m'(t)=-g(y(t))a_m(t)\\
y_{m+1}'(t)=h(y(t))b(t)\,,
\end{array} \right.
\end{equation}
where the state $y:=(y_1^T,y_2^T,...,y_{m+1}^T)^T\in\mathbb{R}^N$ and every $y_i\in\mathbb{R}^n$. We assume  that $g:\mathbb{R}^n\rightarrow \mathbb{R}^+$ and $h:\mathbb{R}^n\rightarrow \mathbb{R}^+$, so every pursuer has the same velocity rule  (depending from the position; in the previous case $y$ was the distance between agents). This will give us a geometric property which will guarantee the decomposability of the problem. The appropriate target set in the present setting is
\begin{multline}
\mathcal{T}=\left\{ (y_1^T,y_2^T,...,y_{m+1}^T)\in \mathbb{R}^N: \right.\\
\left.\min_{i\in \{1,2...m\}}|y_i-y_{m+1}|\leq r \right\}
\end{multline}
Fig. \ref{targset} illustrates the target set in the case $m=2$, $n=1$.\par 

The Hamiltonian, which is
%We can easily show the non convexity of the Hamiltonian: with a proof analogous to Proposition \ref{convex} we can show that 
\begin{equation*}
H(x,p)=g(x) \sum_{i=1}^m  \rho_a |p_i| - h(x)\rho_b |p_{m+1}|-1\,,
\end{equation*}
is not convex in the $p$ variable.

\begin{theorem}\label{decom2}
Assume \eqref{DYN}, \eqref{H3}. Let $u_i(.):\mathbb{R}^n\rightarrow \mathbb{R}$ the solution of the following equation
\begin{equation}
\left\{
\begin{array}{l}
v_i(x_{i})+\max\limits_{a_i}\min\limits_{b}\left\{(g(x_{i})a_i, -h(x)b)\cdot  D v_i(x_{i})\right\}=1 \\
\phantom{ggggggggggggggggggggggg} 
x\in (r,+\infty]^{2n}\\
v(x_{i})=0 
\phantom{gggggggggggggggg} 
x_{i}\in \{|x_1-x_{2}|\leq r \}
\end{array} \right.
\end{equation} 
with $a_i\in B_n(0,\rho_{a})$, $b\in B_n(0,\rho_{b})$. Define the function $u_i(x):\mathbb{R}^N\rightarrow\mathbb{R}$ as
\begin{equation}
u_i(x=(x_1,\ldots,,x_m))=v_i(x_i) \;.
%\quad \forall (x_1,...x_m)\in\mathbb{R}^{N}
\end{equation}
Then the value function for the original problem is
$$ u(x)=\min\{u_1(x),\ldots,u_{m}(x)\}$$
\end{theorem}
\medskip

\proof
To simplify the notation, suppose that $n=1$. (The case for a general $n$ is treated in the same way). We can see that the Hamiltonian is convex in $p$ along rays in every direction with exception of the $e_{m+1}$ direction. (Here, $e_{i}$ is $i$'th canonical basis vector). \par
We can repeat the main steps in the proof of Theorem \ref{decom} with the exception of the verification of condition $(C)$. In this case  we know, for geometric reason,  the superdifferential of the function $u(.)$ is an element aligned with $e_{m+1}$, i.e.
$$ \xi_i, \xi_j\in \partial^F  u(x), \quad (\xi_i-\xi_j)\cdot e_{m+1}=0.$$
This follows from the fact that $e_{m+1}$ is tangential to the switching  interface. Writing $u_{i}$ and $u_j$ for two reduced value functions, we can show that $u_i(.)=u_j(.)$, i.e. they solve the same equation. It follows that the switching interface is located where two reduced value functions 
% in $\mathbb{R}^N$
 coincide, i.e. where $u_j(x)=u_i(x)$. Writing $n$ for the normal of the switching interface, we have 

\begin{multline}
 n\cdot e_{m+1}=\left( 0, \hdots, 0, \frac{\partial u_i}{\partial x_i} , 0 ,\hdots , 0 , \frac{\partial u_i}{\partial x_j} , 0 ,\hdots ,0\right)^t \\
\cdot \left(  0 , \hdots , 0, 1  \right)^t  =0 \;.
\end{multline}

Condition $(C)$ can now be validated. the state representation of $u(.)$ is therefore valid by Thm. 3.1.
\endproof

As an example, consider the $n=2$ case. An evader has position denoted by $x_e\in\mathbb{R}^2$  and $m\in \mathbb{N}$ pursuers have positions denoted by  $x_1, x_2, ... x_n \in \mathbb{R}^2$. Take $m=3$. It is assumed that
% this example we consider 
\begin{equation}
\left\{
\begin{array}{ll}
    g(x)=1-0.5 \cos(\pi x) &\hbox{ if } |x_2|<0.5\\
    g(x)=1 \quad & \hbox{elsewhere}
\end{array} \right.
\end{equation}
and $h(x)\equiv 0.4$. \par

\begin{figure}[th]
\begin{center}
\includegraphics[height=7cm]{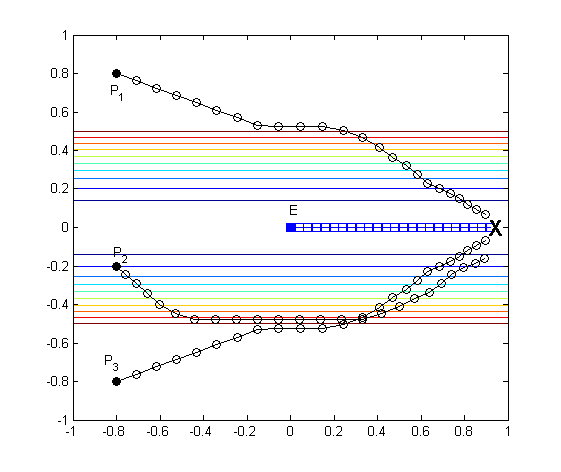}
\caption{Test 2: optimal trjectories for a 3pursuers game.} \label{fig8}
\end{center}
\end{figure}

\begin{figure}[th]
\begin{center}
\includegraphics[height=7cm]{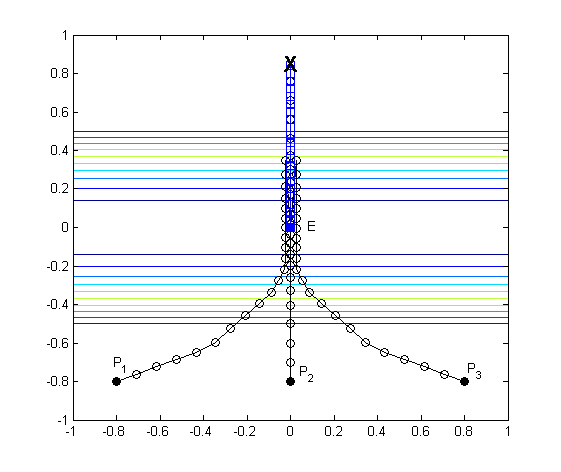}
\caption{Test 2: optimal trjectories for a 3pursuers game.} \label{fig9}
\end{center}
\end{figure}

Figs. \ref{fig8}, \ref{fig9} show some optimal trajectories and the level sets of the velocity function of the pursuers.

%\begin{figure}[th]
%\begin{center}
%\includegraphics[height=6cm]{value.png}
%\caption{Test 2: value function projected in a 3D space ($x_4=0$).} \label{fig8}
%\end{center}
%\end{figure}

\section{ACKNOWLEDGMENTS}

This work was supported by the European Union under the 7th Framework Programme FP7-PEOPLE-2010-ITN SADCO, Sensitivity Analysis for Deterministic Controller Design.

\end{document}